\documentclass[12pt]{amsart}

\usepackage{amssymb}  
\usepackage{latexsym} 
\usepackage{comment}

\usepackage[all]{xy}

\DeclareFontEncoding{OT2}{}{} 
\newcommand{\textcyr}[1]{%
 {\fontencoding{OT2}\fontfamily{cmr}\fontseries{m}\fontshape{n}\selectfont #1}}

\newcommand{\Sha}{{\mbox{\textcyr{Sh}}}}

\newcommand{\Z}{{\mathbb Z}}
\newcommand{\Q}{{\mathbb Q}}
\newcommand{\R}{{\mathbb R}}
\newcommand{\F}{{\mathbb F}}
\newcommand{\BP}{{\mathbb P}}

\newcommand{\CO}{{\mathcal O}}

\newcommand{\CE}{{\mathcal E}}

\newcommand{\To}{\longrightarrow}
\newcommand{\Too}{\;\longrightarrow\;}

\newcommand{\tensor}{\otimes}

\newcommand{\GL}{\operatorname{GL}}
\newcommand{\SL}{\operatorname{SL}}
\newcommand{\PGL}{\operatorname{PGL}}

\newcommand{\Mat}{\operatorname{Mat}}
\newcommand{\Gal}{\operatorname{Gal}}
\newcommand{\Aut}{\operatorname{Aut}}
\newcommand{\End}{\operatorname{End}}
\newcommand{\Hom}{\operatorname{Hom}}
\newcommand{\Map}{\operatorname{Map}}

\newcommand{\Sel}{\operatorname{Sel}}

\newcommand{\Sym}{\operatorname{Sym}}
\newcommand{\WC}{\operatorname{WC}}
\newcommand{\Cl}{\operatorname{Cl}}
\newcommand{\Gm}{\mathbb{G}_m}

\newcommand{\inj}{\hookrightarrow}

\newcommand{\res}{\operatorname{res}}
\newcommand{\unr}{{\text{\rm unr}}}

\newcommand{\diw}{\operatorname{div}}

\newcommand{\rank}{\operatorname{rank}}

\newcommand{\tors}{{\text{\rm tors}}}

\newcommand{\Br}{\operatorname{Br}}
\newcommand{\inv}{\operatorname{inv}}
\newcommand{\Ob}{\operatorname{Ob}}

\newenvironment{Proof}{\par\noindent{\sc Proof:}}%
                      {\hspace*{\fill}\nobreak$\Box$\par}
\newenvironment{Remark}{\par\noindent{\sc Remark:}}{\par}

   {\begin{list}{}{\settowidth{\labelwidth}{#1}%
                   \setlength{\leftmargin}{\labelwidth}%
                   \addtolength{\leftmargin}{\labelsep}}}%
   {\end{list}}

\newtheorem{Theorem}{Theorem}[section]

\newtheorem{Proposition}[Theorem]{Proposition}
\newtheorem{Corollary}[Theorem]{Corollary}

\numberwithin{equation}{section}

\begin{document}

\title{Descent on elliptic curves}

\author{Michael Stoll}
\address{School of Engineering and Science,
         International University Bremen 
         (Jacobs University Bremen as of spring 2007),
         P.O.Box 750561,
	 28725 Bremen, Germany.}
\email{m.stoll@iu-bremen.de}
\date{November 22, 2006}

\begin{abstract}
  Let $E$ be an elliptic curve over~$\Q$ (or, more generally, a number
  field). Then on the one hand, we have the finitely generated abelian
  group $E(\Q)$, on the other hand, there is the Shafarevich-Tate group
  $\Sha(\Q, E)$. {\em Descent} is a general method of getting information
  on both of these objects --- ideally complete information on the 
  Mordell-Weil group $E(\Q)$, and usually partial information on~$\Sha(\Q, E)$.
  
  What descent does is to compute (for a given $n > 1$) the {\em $n$-Selmer
  group} $\Sel^{(n)}(\Q, E)$; it sits in an exact sequence
  \[ 0 \To E(\Q)/nE(\Q) \To \Sel^{(n)}(\Q, E) \To \Sha(\Q, E)[n] \To 0 \]
  and thus contains combined information on $E(\Q)$ and~$\Sha(\Q, E)$.
  
  The main problem I want to discuss in this ``short course'' is how to
  actually do this explicitly, with some emphasis on obtaining representations
  of the elements of the Selmer group as explicit covering spaces of~$E$.
  These explicit representations are useful in two respects --- they allow
  a search for rational points (if successful, this proves that the element
  is in the image of the left hand map above), and they provide the starting
  point for performing ``higher'' descents (e.g., extending a $p$-descent
  computation to a $p^2$-descent computation).
  
  \medskip
  
  {\bf Prerequisites:} Basic knowledge of elliptic curves (e.g., Silverman's
  book~\cite{Silverman}), some Galois cohomology and algebraic number theory 
  (e.g., Cassels--Fr\"ohlich~\cite{CasselsFroehlich}).
\end{abstract}

\maketitle

\renewcommand{\baselinestretch}{1.1}
\renewcommand{\arraystretch}{1.3}

\pagebreak


The results described in these notes (if not ``classical'', i.e., to
be found in, e.g., Silverman's book~\cite{Silverman}) 
were obtained in collaboration
with John Cremona, Tom Fisher, Cathy O'Neil, Ed Schaefer and Denis Simon.
See the papers \cite{SchaeferStoll,Paper1,Paper2,Paper3} for a more detailed account.

\section{The Selmer Group}

In the following, $K$ will be a number field, and $E$ will be an elliptic
curve defined over~$K$. $E$ is an algebraic group over~$K$, and so its
set of rational points, $E(K)$, forms a group,
the so-called Mordell-Weil group. By the Mordell-Weil theorem,
it is a finitely
generated abelian group, and one of the big questions is how to determine it
(in the sense of, say, giving generators as points in~$E(K)$ and relations).
Descent is the main tool used for that, both in theory and in practice.
Doing an $n$-descent on~$E$ means to compute the $n$-Selmer group
$\Sel^{(n)}(K,E)$, which we will introduce in this section.

Note that saying that $E(K)$ is a finitely generated abelian group
amounts to asserting the existence of an exact sequence
\[ 0 \Too E(K)_\tors \Too E(K) \Too \Z^r \Too 0 \]
with $r \ge 0$ an integer and $E(K)_\tors$ a finite abelian group;
it consists of all elements of~$E(K)$ of finite order.
Less canonical, but sometimes more convenient, we also have
\[ E(K) \cong E(K)_\tors \oplus \Z^r \,. \]

For any concrete curve~$E$, it is fairly straightforward to find~$E(K)_\tors$,
and we will not be concerned with how to do that in these lectures.
The hard part is to determine the rank~$r$. This is where descent helps.


\subsection{Definition and first properties}\strut

Let $n > 1$ be an integer. The usual definition of the $n$-Selmer group
makes use of Galois cohomology. Consider the short exact sequences
of $G_K = \Gal(\bar{K}/K)$-modules
\[ 0 \Too E[n](\bar{K}) \Too E(\bar{K}) \stackrel{n}{\Too} E(\bar{K}) \Too 0 \]
(which is usually just written
\[ 0 \Too E[n] \Too E \stackrel{n}{\Too} E \Too 0 \,). \]
Then we have the long exact sequence of cohomology groups
\[ 0 \To E[n](K) \To E(K) \stackrel{n}{\To} E(K)
    \stackrel{\delta}{\To} H^1(K, E[n]) \To H^1(K, E) 
    \stackrel{n}{\To} H^1(K, E) \,.
\]
We deduce from it another short exact sequence:
\[ 0 \Too E(K)/n E(K) \stackrel{\delta}{\Too} H^1(K, E[n])
     \Too H^1(K, E)[n] \Too 0 
\]

It turns out that knowing $E(K)$ is essentially equivalent to knowing
its free abelian rank~$r = \rank E(K)$. (Once we know~$r$, we can look
for points until we have found $r$ independent ones. Then we only need
to find the $K$-rational torsion points and ``saturate'' the subgroup
generated by the independent points. All of this can be done effectively.)
Now the idea is to use the above exact sequence to at least get an
upper bound on~$r$: $r$ can be read off from the size of the group
$E(K)/n E(K)$ on the left, and so any bound on that group will provide
us with a bound on~$r$. From the exact sequence, we see that 
$E(K)/n E(K)$ sits inside $H^1(K, E[n])$; however this group is infinite,
and so it does not give a bound.

But we can use some additional information. We know (trivially) that
any $K$-rational point on~$E$ is also a $K_v$-rational point, for all
places $v$ of~$K$. Now it is possible to compute the image of the
local map
\[ E(K_v)/n E(K_v) \stackrel{\delta_v}{\Too} H^1(K_v, E[n]) \]
for any given~$v$ explicitly; and for all but a finite explicitly
determinable set of places~$S$, the image just consists of the
``unramified part'' of~$H^1(K_v, E[n])$. This means that in some sense, we
can compute all the necessary ``local'' conditions and use this
information in bounding the ``global'' group $E(K)/n E(K)$.
Formally, we define the {\em $n$-Selmer group} of~$E$, $\Sel^{(n)}(K, E)$,
to be the subgroup of $H^1(K, E[n])$ of elements that under all
restriction maps $\res_v$ are in the image of~$\delta_v$ in the
following diagram.
\[ \SelectTips{cm}{}
   \xymatrix{ 0 \ar[r] & 
              E(K)/n E(K) \ar[dd] \ar[r]^{\delta} &
              H^1(K, E[n]) \ar[dd]^{\prod_v \res_v} \ar[ddr]^{\alpha} \ar[r] &
              H^1(K, E)[n] \ar[dd]^{\prod_v \res_v} \ar[r] &
              0 \\
              \\
              0 \ar[r] &
              \prod\limits_v E(K_v)/n E(K_v) \ar[r]^{\prod_v \delta_v} & 
              \prod\limits_v H^1(K_v, E[n]) \ar[r] &
              \prod\limits_v H^1(K_v, E)[n] \ar[r] &
              0 }
\]
Equivalently, $\Sel^{(n)}(K, E)$ is the kernel of the map~$\alpha$.
The image of~$\Sel^{(n)}(K, E)$ in~$H^1(K, E)[n]$ is the kernel of
the rightmost vertical map in the diagram. More generally, one defines
the {\em Shafarevich-Tate group} of~$E$, $\Sha(K, E)$ to be
\[ \Sha(K, E) = \ker\bigl(H^1(K, E) \To \prod_v H^1(K_v, E)\bigr) \,. \]
Then we get another short exact sequence:
\[ 0 \Too E(K)/n E(K) \stackrel{\delta}{\Too}
     \Sel^{(n)}(K, E) \Too \Sha(K, E)[n] \Too 0 \,.
\]
This time, one can (and we will) prove that the middle group is finite.
And at least in principle, it is computable. In this way, we can compute
the product $(\# E(K)/n E(K))(\# \Sha(K, E)[n])$, and in particular,
we obtain a bound on the rank~$r$. The obstruction against this bound
being sharp lies in $\Sha(K, E)$, which is therefore also an interesting
object. Of course, its size (conjectured, but not generally proved to be
finite) also shows up in the famous Birch and Swinnerton-Dyer conjecture,
and there are other reasons to study $\Sha(K, E)$ for its own sake.

We need some more notions and notation. The {\em unramified part}
of $H^1(K_v, E[n])$ is the kernel of the restriction map
$H^1(K_v, E[n]) \To H^1(K_v^\unr, E[n])$. For any finite set of places~$S$
of~$K$ containing the infinite places, we define $H^1(K, E[n]; S)$ to
be the subgroup of~$H^1(K, E[n])$ of elements that map into the unramified
part of~$H^1(K_v, E[n])$ for all places $v \notin S$.

The finiteness of the Selmer group then follows from the two observations
that $\Sel^{(n)}(K, E) \subset H^1(K, E[n]; S)$ for a suitable finite set~$S$,
and that $H^1(K, E[n]; S)$ is finite for all finite sets~$S$ of places of~$K$.

The latter is a standard fact; in the end it reduces to the two basic
finiteness results of algebraic number theory: finiteness of the class group
and finite generation of the unit group.

\begin{Theorem}
  If $S$ is a finite set of places
  of~$K$ containing the infinite places, then $H^1(K, E[n]; S)$ is finite.
\end{Theorem}

\begin{Proof}
  There is a finite extension $L = K(E[n])$ of~$K$ (the $n$-division
  field of~$E$) such that $E[n]$ becomes a trivial $L$-Galois module.
  We have the inflation-restriction exact sequence
  \[ 0 \Too H^1(L/K, E[n](L)) \Too H^1(K, E[n]) \Too H^1(L, E[n]) \,, \]
  and the group on the left is finite. Taking into account the ramification
  conditions, we see that $H^1(K, E[n]; S)$ maps into $H^1(L, E[n]; S_L)$
  with finite kernel, where $S_L$ is the set of places of~$L$ above some
  place of~$K$ in~$S$. Therefore it suffices to show that $H^1(L, E[n]; S_L)$
  is finite. Now 
  \[ H^1(L, E[n]) = H^1(L, (\Z/n\Z)^2) = \Hom(G_L, (\Z/n\Z)^2) \,, \]
  and the ramification condition means that the fixed field of the kernel
  of a homomorphism coming from $H^1(L, E[n]; S_L)$ is unramified 
  outside~$S_L$. On the other hand, this fixed field is an abelian extension
  of exponent dividing~$n$; it is therefore contained in the maximal
  abelian extension~$M$ of exponent~$n$ that is unramified outside~$S_L$.
  
  By Kummer theory ($L$ contains the $n$th roots of unity because of the
  $n$-Weil pairing), $M = L(\sqrt[n]{U})$ for some subgroup 
  $U \subset L^\times/(L^\times)^n$. Enlarging $S_L$ by including the
  primes dividing~$n$, the ramification condition translates into
  \[ U = L(S_L, n)
       = \{\alpha \in L^\times
             : n \mid v(\alpha) \text{\ for all $v \notin S_L$}\}/(L^\times)^n
  \]
  (the ``$n$-Selmer group of $\CO_{L,S_L}$''). Applying the Snake Lemma
  to the diagram below then provides us with the exact sequence
  \[ 0 \Too \CO_{L,S_L}^\times/(\CO_{L,S_L}^\times)^n
       \Too L(S_L, n) \Too \Cl_{S_L}(L)[n] \Too 0 \,.
  \]
  Since the $S_L$-unit group $\CO_{L,S_L}^\times$ is finitely generated
  and the $S_L$-class group $\Cl_{S_L}(L)$ is finite, $U = L(S_L, n)$
  is finite, and hence so is the extension~$M$. We see that all the 
  homomorphisms have to factor through the finite group $\Gal(M/L)$,
  whence $H^1(L, E[n]; S_L)$ is finite.
  \[ \SelectTips{cm}{}
     \xymatrix{ & 0 \ar[d] & 0 \ar[d] & 0 \ar[d] & \\
                & \CO_{L,S_L}^\times \ar[r]^n \ar[d] &
                  \CO_{L,S_L}^\times \ar[r] \ar[d] &
                  L(S_L, n) \ar[d] \ar@{-->} [r] & \\
                & L^\times \ar[r]^n \ar[d] & L^\times \ar[r] \ar[d] &
                  L^\times/(L^\times)^n \ar[d] \ar[r] & 0 \\
               0 \ar[r] & I_{S_L} \ar[r]^n \ar[d] &
                  I_{S_L} \ar[r] \ar[d] & I_{S_L}/n I_{S_L} \ar[r] & 0 \\
               \ar@{-->} [r] & \Cl_{S_L}(L) \ar[r]^n \ar[d] &
                \Cl_{S_L}(L) \ar[d] & & \\
                & 0 & 0 & &
              }
  \]
\end{Proof}

The next result implies that $\Sel^{(n)}(K, E)$ is finite.

\begin{Theorem} \label{Ramification}
  $\Sel^{(n)}(K, E) \subset H^1(K, E[n]; S)$ where $S$ is any finite set
  of places of~$K$ containing the infinite places, the places dividing~$n$
  and the finite places $v$ such that $\gcd(c_v(E), n) > 1$, where
  $c_v(E)$ is the Tamagawa number of~$E$ at~$v$.
\end{Theorem}

\begin{Proof}
  We have to show that for $\xi \in \Sel^{(n)}(K, E) \subset H^1(K, E[n])$
  and for $v \notin S$, $\xi$ maps to zero in $H^1(K_v^\unr, E[n])$.
  Consider the exact sequences
  \begin{gather*}
    0 \Too E(K_v^\unr)^0 \Too E(K_v^\unr) \Too \Phi_v(\bar{k}_v) \Too 0 \\
    0 \Too E(K_v^\unr)^1 \Too E(K_v^\unr)^0 \Too \CE(\bar{k}_v)^0 \Too 0
  \end{gather*}
  Here $\CE$ is the N\'eron model of~$E$ over $\CO_{K_v}$, $\CE(\bar{k}_v)^0$
  is the connected component of the identity on the special fiber of~$\CE$,
  $E(K_v^\unr)^0$ is the subgroup of points mapping into $\CE(\bar{k}_v)^0$
  (the points of good reduction on a minimal Weierstrass model at~$v$), and
  $E(K_v^\unr)^1$ is the kernel of reduction at~$v$. Applying the Snake
  Lemma to multiplication-by-$n$ on these sequences gives exact sequences
  of cokernels
  \[ E(K_v^\unr)^0/n E(K_v^\unr)^0
       \Too E(K_v^\unr)/n E(K_v^\unr) \Too \Phi_v(\bar{k}_v)/n \Phi_v(\bar{k}_v)
  \]
  and
  \[ E(K_v^\unr)^1/n E(K_v^\unr)^1 \Too E(K_v^\unr)^0/n E(K_v^\unr)^0
      \Too \CE(\bar{k}_v)^0/n \CE(\bar{k}_v)^0 \,.
  \]
  We claim that the image of $E(K_v)$ in~$E(K_v^\unr)$ is divisible by~$n$.
  Let $P \in E(K_v)$. Then in the first sequence, the image of~$P$ in the
  group furthest on the right is in 
  $\Phi_v(k_v)/(\Phi_v(k_v) \cap n \Phi_v(\bar{k}_v))$,
  and this group is trivial, since $c_v(E) = \#\Phi_v(k_v)$ is prime to~$n$.
  Hence the image of~$P$ comes from $E(K_v^\unr)^0/n E(K_v^\unr)^0$.
  In the second sequence, the group on the right is trivial, because the
  $\bar{k}_v$-points of an algebraic group over~$k_v$ form a divisible group.
  The first group is also trivial, because it is a $\Z_p$-module
  (with $p$ the residue characteristic of~$v$),
  and $n$ is invertible in this ring. Hence $E(K_v^\unr)^0/n E(K_v^\unr)^0 = 0$,
  and the image of $P$ must vanish in~$E(K_v^\unr)/n E(K_v^\unr)$.
  Now consider the following diagram.
  \[ \SelectTips{cm}{}
     \xymatrix{ E(K_v)/n E(K_v) \ar[r]^-{\delta_v} \ar[d]
                  & H^1(K_v, E[n]) \ar[d] \\
                E(K_v^\unr)/n E(K_v^\unr) \ar[r]^-{\delta_v^\unr} 
                  & H^1(K_v^\unr, E[n])
              }
  \]
  We have seen that the left vertical arrow is the zero map, therefore the
  image of $\delta_v$ also maps trivially under the right vertical map.
  This exactly means that the elements of the Selmer group (mapping into
  the image of~$E(K_v)$ in~$H^1(K_v, E[n])$) are unramified at~$v$.
\end{Proof}

\begin{Remark}
  The proof shows that in general, the image of~$E(K_v)$ 
  in~$H^1(K_v^\unr, E[n])$ is isomorphic to 
  $\Phi_v(k_v)/(\Phi_v(k_v) \cap n\Phi_v(\bar{k}_v))$ for all finite
  places~$v$ that do not divide~$n$. In particular, the order of the
  image divides the Tamagawa number~$c_v(E)$.
\end{Remark}


\subsection{Interpretation of Selmer group elements}\strut

The definition of the Selmer group as a subgroup of~$H^1(K, E[n])$ is
rather abstract, its elements being given by classes of 1-cocycles with
values in~$E[n]$. However, it is possible to give the Selmer group elements
much more concrete interpretations. This is based on the following general
fact.

\begin{Proposition}
  Let $X$ be some sort of algebraic or geometric object, defined over~$K$.
  Then the set of {\em twists} of~$X$, i.e., objects $Y$ defined over~$K$
  such that $X$ and~$Y$ are isomorphic over~$\bar{K}$, up to $K$-isomorphism,
  is parametrized by $H^1(K, \Aut_{\bar{K}}(X))$. (When the automorphism
  group of~$X$ is abelian, this is an abelian group; otherwise, it is
  a pointed set with the class of~$X$ as its distinguished element.)
\end{Proposition}

\begin{Proof}
  This is quite standard (at least in many concrete manifestations).
  The map from the twists to $H^1(K, \Aut_{\bar{K}}(X))$ is obtained as
  follows. Let $Y$ be a twist of~$X$. Then there is an isomorphism
  $\phi : Y \to X$, defined over~$\bar{K}$. Then 
  $\xi_\sigma = \phi^\sigma \phi^{-1}$ defines a 1-cocycle with values
  in~$\Aut_{\bar{K}}(X)$, and postcomposing $\phi$ by an automorphism
  of~$X$ changes $\xi$ into a cohomologous cocycle.
  To get the map in the reverse direction, one takes the $\bar{K}$-``points''
  of~$X$ and ``twists'' the action of $G_K$ by~$\xi$ by decreeing that
  the action of~$\sigma$ on $Y$ (which has the same underlying set of
  $\bar{K}$-``points'' as~$X$) is given by the action of~$\sigma$ on~$X$,
  followed by~$\xi_{\sigma^{-1}}$.
\end{Proof}

Since $\Sel^{(n)}(K, E) \subset H^1(K, E[n])$, this means that we can obtain
interpretations of Selmer group elements via interpretations of elements
of $H^1(K, E[n])$ as twists. So we have to look for ``objects'' whose
($\bar{K}$-)automorphism group is~$E[n]$.


\subsection*{Principal homogeneous spaces}

Let us first look at a somewhat simpler situation related to~$H^1(K, E)$.
Consider ``objects'' of the form 
$\SelectTips{cm}{}\xymatrix@1{C \ar@{-->}[r]^\cong & E}$, where the
isomorphism is defined over~$\bar{K}$. Two such diagrams are isomorphic
if there is an isomorphism $C \stackrel{\cong}{\To} C'$ and a point $P \in E$
such that the diagram
\[ \SelectTips{cm}{}
   \xymatrix{ C \ar[d] \ar@{-->}[r] & E \ar@{-->}[d]^{\cdot\; + P} \\
              C' \ar@{-->}[r] & E
            }
\]
commutes. Then the automorphisms of the trivial object $E \stackrel{=}{\To} E$
are just the translations, so the automorphism group is~$E(\bar{K})$.
The objects are called {\em principal homogeneous spaces} for~$E$, and they
are classified (up to $K$-isomorphisms) by the {\em Weil-Ch\^atelet group} 
$\WC(K, E) = H^1(K, E)$. 

Given a curve~$C$ as above, we can change the isomorphism to~$E$ by any
translation without changing the isomorphism class of $C$ as a principal
homogeneous space. So given $C$, the only ambiguity in endowing it with
a structure as a principal homogeneous space comes from the automorphism
group of~$E$ as an elliptic curve. Generically, this is just $\{\pm 1\}$,
and so there will be at most two structures as a principal homogeneous
space for~$E$ on~$C$. (The two will coincide when either one has order
dividing~$2$ in~$H^1(K, E)$.)

Note also that a principal homogeneous space has an algebraic group action 
of~$E$ on it. If $\phi : C \To E$ is an isomorphism (over~$\bar{K}$), then 
\[ C \times E \ni (P, Q) \longmapsto P + Q := \phi^{-1}(\phi(P) + Q) \in C \]
is defined over~$K$, since it is unchanged when $\phi$ is post-composed
with a translation on~$E$. Also, there is a well-defined (over~$K$)
``difference morphism''
\[ C \times C \ni (P, P') \longmapsto P - P' := \phi(P) - \phi(P') \in E \]
such that, for example, $P + (P' - P) = P'$. Conversely, given morphisms
$C \times E \To C$ and $C \times C \To E$ satisfying the usual properties,
$C$ becomes a principal homogeneous space (in the sense above) by picking
any point $P_0 \in C$ and considering the isomorphism 
$C \ni P \mapsto P - P_0 \in E$. So we could have defined principal 
homogeneous spaces also through actions of~$E$ on curves~$C$. (In fact,
this is what is usually done.)

If $C$ has a $K$-rational point~$P$, then there is a $K$-defined
isomorphism $\phi : C \To E$ that maps $P$ to~$O$. We obtain a diagram
\[ \SelectTips{cm}{}
   \xymatrix{ C \ar[d]_\phi \ar[r]^\phi & E \ar@{=}[d] \\
              E \ar@{=}[r] & E
            }
\]
showing that $C$ is trivial as a principal homogeneous space (i.e.,
$K$-isomorphic to $E \stackrel{=}{\To} E$).

In this context, the elements of $\Sha(K, E)$ are represented by
principal homogeneous spaces with $K_v$-points for every place~$v$,
or with points ``everywhere locally'', up to 
isomorphism over~$K$. Nontrivial elements 
of~$\Sha(K, E)$ are those that have points everywhere locally, but
no global points, i.e., those that ``fail the Hasse Principle''.


\subsection*{First interpretation: $n$-coverings}

Here, our object $X$ is the multiplication-by-$n$ map $E \stackrel{n}{\To} E$.
The twists are covering maps $C \stackrel{\pi}{\To} E$ such that there is an 
isomorphism $C \To E$ over~$\bar{K}$ such that the following diagram commutes.
\[ \SelectTips{cm}{}
   \xymatrix{ C \ar[r]^\pi \ar@{-->}[d] & E \ar@{=}[d] \\
              E \ar[r]^n & E
            }
\]
Such a $C \stackrel{\pi}{\To} E$ is called an {\em $n$-covering} of~$E$.
An isomorphism between two $n$-coverings $C \stackrel{\pi}{\To} E$
and $C' \stackrel{\pi'}{\To} E$ is given by an isomorphim 
$\phi : C \To C'$ such that the following diagram commutes.
\[ \SelectTips{cm}{}
   \xymatrix{ C \ar[r]^\pi \ar[d]_\phi & E \ar@{=}[d] \\
              C' \ar[r]^{\pi'} & E
            }
\]
The automorphisms of $E \stackrel{n}{\To} E$ are then given by the
translations by $n$-torsion points (acting on the left~$E$), so that we
indeed obtain $E[n]$ as the automorphism group.

In this interpretation, the map $E(K)/n E(K) \To H^1(K, E[n])$ comes about
as follows. To a point $P \in E(K)$, we associate the $n$-covering
$E \stackrel{\pi}{\To} E$ such that $\pi(Q) = nQ + P$. It is easy to check
that the isomorphism class of the covering only depends on $P$ mod~$n E(K)$.
On the other hand, each $n$-covering $C \stackrel{\pi}{\To} E$ such that
$C(K) \neq \emptyset$ is isomorphic to one of this form: there is an
isomorphism between $C$ and~$E$ defined over~$K$ (mapping
a $K$-rational point on~$C$ to $O \in E$); under the composed map
$E \stackrel{\cong}{\To} C \stackrel{\pi}{\To} E$, the origin $O$ maps
to some $P \in E(K)$, and then the map must be $Q \mapsto nQ + P$.
Tracing through the definition of the connecting map~$\delta$ shows
that the map defined here coincides with~$\delta$ under our interpretation.

For the Selmer group elements, this means that they correspond to the
$n$-coverings that have points {\em everywhere locally}.

Note that the curve $C$ in an $n$-covering $C \stackrel{\pi}{\To} E$
carries the structure of a {\em principal homogeneous space} for~$E$:
any $\bar{K}$-isomorphism $C \To E$ in the definition above provides
such a structure, and since these isomorphisms are all related by
translations (by $n$-torsion points), the isomorphism class of the
principal homogeneous space structure is well-defined. This gives us
the map $H^1(K, E[n]) \To H^1(K, E)$ in the coverings interpretation.


\subsection*{Second interpretation: Maps to $\BP^{n-1}$}

Here is another interpretation. On~$E$, we can consider the map to~$\BP^{n-1}$
that is given by the complete linear system associated to~$n \cdot O$.
Other objects are diagrams
\[ \SelectTips{cm}{}
   \xymatrix{ C \ar[r] \ar@{-->}[d]_{\cong} & S \ar@{-->}[d]^{\cong} \\
              E \ar[r] & \BP^{n-1}
            }
\]
with the dashed isomorphisms defined over~$\bar{K}$. Twists of~$\BP^{n-1}$
like~$S$ are called {\em Severi-Brauer varieties}; they are classified
(according to the general principle) by $H^1(K, \PGL_n)$. Since $\PGL_n$
is non-abelian, this is just a pointed set; however, applying Galois
cohomology to the exact sequence
\[ 0 \Too \Gm \Too \GL_n \Too \PGL_n \Too 0 \,, \]
one obtains an injection $H^1(K, \PGL_n) \To H^2(K, \Gm) = \Br(K)$
identifying $H^1(K, \PGL_n)$ with the $n$-torsion $\Br(K)[n]$ in the
Brauer group of~$K$.

An isomorphism between two such diagrams is given by a pair of 
isomorphisms $C \To C'$ and $S \To S'$ such that the diagram
\[ \SelectTips{cm}{}
   \xymatrix{ C \ar[rrr] \ar[ddd] \ar@{-->}[dr] & & & 
                     S \ar[ddd] \ar@{-->}[dl] \\
               & E \ar[r] \ar@{-->}[d]^{\cdot\;+P} & \BP^{n-1} \ar@{-->}[d] \\
               & E \ar[r] & \BP^{n-1} \\
              C' \ar[rrr] \ar@{-->}[ur] & & & S' \ar@{-->}[ul]
            }
\]
commutes, with some choice of $P \in E$ and some automorphism of~$\BP^{n-1}$ 
over~$\bar{K}$. Automorphisms of $E \To \BP^{n-1}$ are therefore given by 
translations on~$E$ that fix the linear system $|n \cdot O|$. Translation
by~$P$ does this if and only if $n \cdot P \in |n \cdot O|$, i.e., iff
$P \in E[n]$. So we obtain again the correct automorphism group, and we
see that $H^1(K, E[n])$ also classifies diagrams as above, up to isomorphism
over~$K$. The map to $H^1(K, E)$ comes through ``forgetting'' the right
half of the diagram.

In particular, we see that the curve $C$ in both the coverings and the
maps to $\BP^{n-1}$ interpretation of a given element in~$H^1(K, E[n])$
is the same (as a principal homogeneous space for~$E$). So, also in
our second interpretation, the elements of the Selmer group are those
diagrams such that $C$ has points everywhere locally. This implies that
also the Severi-Brauer variety~$S$ has points everywhere locally.
Now there is the very important local-global principle for the Brauer
group:
\[ 0 \Too \Br(K) \Too \bigoplus_v \Br(K_v) 
     \stackrel{\sum_v \inv_v}{\Too} \Q/\Z \Too 0
\]
is exact. In particular, an element of the Brauer group of~$K$ that is locally
trivial is already (globally) trivial. This implies that $S$ has a
$K$-rational point, and so $S \cong \BP^{n-1}$. Whence the following result.

\begin{Proposition}
  The elements of~$\Sel^{(n)}(K, E)$ are in 1-to-1 correspondence with
  $K$-isomorphism classes of diagrams
  \[ \SelectTips{cm}{}
     \xymatrix{ C \ar[r] \ar@{-->}[d]_{\cong}
                  & \BP^{n-1} \ar@{-->}[d]^{\cong} \\
                E \ar[r] & \BP^{n-1}
              }
  \]
  such that $C$ has points everywhere locally.
\end{Proposition}

Let us look at what this means for various small values of~$n$.
\begin{itemize}\addtolength{\itemsep}{1mm}
  \item[$n = 2$:]
    On~$E$, the map to~$\BP^1$ given by~$|2 \cdot O|$ is just the 
    $x$-ccordinate.
    It is a 2-to-1 map ramified in four points (namely, $E[2]$). 
    Any twist $C \To \BP^1$ will have the same geometric properties, which
    means that $C$ can be realized by a model of the form
    \[ y^2 = f(x) = a\,x^4 + b\,x^3 + c\,x^2 + d\,x + e \,. \]
    (This is an affine model; a somewhat better way is to consider
    \[ y^2 = F(x, z) = a\,x^4 + b\,x^3 z + c\,x^2 z^2 + d\,x z^3 + e\,z^4 \]
    in a $(1,2,1)$-weighted projective plane.)
  \item[$n = 3$:]
    The map $E \To \BP^2$ given by~$|3 \cdot O|$ is an embedding of degree~3,
    realizing $E$ as a plane cubic curve. Similarly, any element of the
    3-Selmer group can be realized as a plane cubic curve, with an action
    of~$E[3]$ on it through linear automorphisms.
  \item[$n = 4$:]
    Here we obtain a degree-4 embedding into~$\BP^3$; its image is given
    as the intersection of two quadrics.
\end{itemize}

More generally, for $n \ge 4$, the image in~$\BP^{n-1}$ is given as the
intersection of \hbox{$n(n-3)/2$} quadrics; for $n \ge 5$, this is no longer
a complete intersection. For $n = 5$, there is a nice description by
sub-Pfaffians of a $5 \times 5$ matrix of linear forms.

For general elements of~$H^1(K, E[n])$, we obtain another forgetful map.
If we forget the left hand side of the diagram, then we obtain a map
\[ \Ob : H^1(K, E[n]) \To H^1(K, \PGL_n) = \Br(K)[n] \,, \] 
the ``obstruction''
against an embedding (or a map) into~$\BP^{n-1}$. 

{\bf Warning.} This map is {\em not} a homomorphism!


\subsection{What the Selmer group can be used for}\strut

The most obvious use of the Selmer group is to provide an upper bound
for the Mordell-Weil rank~$r$: we have
\[ n^r = \frac{\#\Sel^{(n)}(K, E)}{\#(E(K)_\tors/n E(K)_\tors) \#\Sha(K, E)[n]}
       \le \frac{\#\Sel^{(n)}(K, E)}{\#(E(K)_\tors/n E(K)_\tors)} \,. \]
To get this, it is sufficient to just compute the order of the Selmer group.

Moreover, by computing the sizes of various Selmer groups (for coprime
values of~$n$), we can compare the bounds we get, and in some cases 
deduce lower bounds on the order of~$\Sha(K, E)$. For example, if we
get $r \le 3$ from the $2$-Selmer group and $r \le 1$ from the $3$-Selmer
group, then we know that $\#\Sha(K, E)[2] \ge 4$.

On the other hand, in order to show that the rank bound we get is sharp,
we need to prove that all elements of the Selmer group come from $K$-rational
points on~$E$. This is rather easy if we find suffiently many independent
points on~$E$. However, in many cases, some of the generators of~$E(K)$
can be rather large and will not be found by a systematic search.
Here, it is useful to represent the elements of the Selmer group as
$n$-coverings $C$. We then have a diagram
\[ \SelectTips{cm}{}
   \xymatrix{ C \ar[rr]^{\deg n} \ar[d]^\pi_{\deg n^2} & & 
                \BP^{n-1} \ar@{.>}[d]^{\deg 2n} \\
              E \ar[rr]^{\deg 2}_x & & \BP^1
            }
\]
with a rational map of degree~$2n$ on the right hand side. 

For example, when $n = 2$, the map $\BP^1 \To \BP^1$ on the right hand
side is given by two quartic forms (the quartic showing up in the
equation of the 2-covering~$C$ and its quartic covariant).

From the general
theory of heights, we expect the logarithmic height of a point on~$C$,
as given by its image in~$\BP^{n-1}$, to be smaller by a factor of about
$1/2n$ than that of the $x$-coordinate of its image on~$E$. This will make
these points much easier to find on~$C$ (in~$\BP^{n-1}$) than on~$E$.
If we find a $K$-rational point on~$C$, we know that the corresponding
element of the $n$-Selmer group is in the image of~$\delta$, and we can
improve the lower bound on the rank. Note that in practice, to really
be fairly certain that the points on~$C$ are as small as expected, it is
necessary to have a ``small'' model of~$C$, i.e., given by equations with
small coefficients.

In case we do not find a $K$-rational point on~$C$, we can use the curve~$C$
as the basis for ``higher descents''; in this way we may be able to prove
that $C$ does not have any $K$-rational points, or find points on curves
that cover~$C$.

The program for the following will therefore be to first show how one
can compute the $p$-Selmer group for a prime number~$p$ (the most important
case). Then we will discuss how to obtain from this computation actual
covering curves. But first, we will introduce another interpretation
of the elements of~$H^1(K, E[n])$.


\subsection*{Third interpretation: Theta groups}

In the second interpretation, on the morphism $E \To \BP^{n-1}$ there is an
action of~$E[n]$; in particular, $E[n]$ acts on~$\BP^{n-1}$ by automorphisms.
Thus we obtain a homomorphism $\chi_E : E[n] \To \PGL_n$. We can then define
$\Theta_E$ by the following diagram.
\[ \SelectTips{cm}{}
   \xymatrix{ 0 \ar[r] & \Gm \ar@{=}[d] \ar[r]^\alpha 
                & \Theta_E \ar[d] \ar[r]^\beta
                & E[n] \ar[d]^{\chi_E} \ar[r] & 0 \\
              0 \ar[r] & \Gm \ar[r] & \GL_n \ar[r] & \PGL_n \ar[r] & 0
            }
\]

\begin{Proposition}
  In the above, for any $\theta, \theta' \in \Theta_E$, we have
  \[ [\theta, \theta'] = \theta \theta' \theta^{-1} {\theta'}^{-1}
       = \alpha\bigl(e_n(\beta(\theta), \beta(\theta')\bigr) \,,
  \]
  where $e_n : E[n] \times E[n] \To \mu_n \inj \Gm$ is the $n$-Weil pairing.
\end{Proposition}

\begin{Proof}
  Let $T, T' \in E[n]$. We have to show that for any two lifts
  $\theta, \theta' \in \GL_n$ of $\chi_E(T)$ and~$\chi_E(T')$, we have
  \[ [\theta, \theta'] = e_n(T, T') \, I_n \]
  (where $I_n$ is the $n \times n$ identity matrix).
  
  For this, note that $\BP^{n-1}$ can be identified with $\BP(L(n \cdot O)^*)$.
  For every $T \in E[n]$, choose $f_T \in \bar{K}(E)^\times$ such that
  $\diw(f_T) = n \cdot T - n \cdot O$. Then the action of~$T$ 
  on~$\BP(L(n \cdot O)^*)$ is induced by 
  \[ L(n \cdot O) \ni f
      \longmapsto \bigl(P \mapsto f_T(P) f(P - T)\bigr) \in L(n \cdot O) \,. 
  \]
  Note that a choice of $\theta$ lifting $\chi_E(T)$ corresponds to a
  choice of~$f_T$. Now the action of the commutator
  $[\theta, \theta']$ is given on~$L(n \cdot O)$ by
  \[ f \longmapsto 
      \bigl(P \mapsto \frac{f_T(P)\,f_{T'}(P - T)}%
                           {f_T(P - T')\,f_{T'}(P)} f(P)\bigr) \,.
  \]
  The factor in front of~$f(P)$ is constant (where defined) and by a standard
  result equal to~$e_n(T, T')$.
\end{Proof}

The proof shows that $\Theta_E$ can be represented as the set
\[ \{(T, f_T) : T \in E[n], f_T \in \bar{K}(E)^\times, 
                \diw(f_T) = n \cdot T - n \cdot O\}
\]
with the group structure given by
\[ (T, f_T) (T', f_{T'}) = (T+T', P \mapsto f_T(P) f_{T'}(P - T)) \,; \]
also $\alpha(\lambda) = (O, \lambda)$ and $\beta(T, f_T) = T$.

More generally, we define a {\em theta group} of level~$n$ for~$E$ to
be an exact sequence (of $K$-group schemes)
\[ 0 \Too \Gm \stackrel{\alpha}{\Too} \Theta
     \stackrel{\beta}{\Too} E[n] \Too 0
\]
such that for $\theta, \theta' \in \Theta$, we have again
\[ [\theta, \theta'] = \alpha\bigl(e_n(\beta(\theta), \beta(\theta')\bigr) \,.
\]
(Note that this implies that $\Theta$ is a {\em central} extension 
of~$E[n]$ by~$\Gm$.)
An isomorphism of two such diagrams is given by a $G_K$-isomorphism
$\phi : \Theta \To \Theta'$ making the following diagram commutative.
\[ \SelectTips{cm}{}
   \xymatrix{ 0 \ar[r] & \Gm \ar@{=}[d] \ar[r]
               & \Theta \ar[d]^\phi \ar[r] & E[n] \ar@{=}[d] \ar[r] & 0 \\
              0 \ar[r] & \Gm \ar[r] & \Theta' \ar[r] & E[n] \ar[r] & 0
            }
\]
Working out what the automorphisms of~$\Theta_E$ are, we find that they
are of the form $(T, f_T) \longmapsto (T, \varphi(T) f_T)$ for a 
homomorphism $\varphi : E[n] \To \Gm$. Such a~$\varphi$ necessarily
takes values in~$\mu_n$, and by the non-degeneracy of the Weil parining~$e_n$,
there is some $T' \in E[n]$ such that $\varphi(T) = e_n(T', T)$.
We see that the automorphism group is again~$E[n]$. Furthermore, we have
the following result.

\begin{Proposition}
  All theta groups of level~$n$ for~$E$ are isomorphic over~$\bar{K}$
  (i.e., as abstract group extensions).
\end{Proposition}

\begin{Proof}
  Short proof: On the level of abstract groups, theta groups are
  classified by $H^2((\Z/n\Z)^2, \bar{K}^\times)$. There is a canonical
  map
  \[ H^2((\Z/n\Z)^2, \bar{K}^\times)
      \To \mathop{\text{\small$\bigwedge$}}\nolimits^{\!2} 
        \Hom((\Z/n\Z)^2, \bar{K}^\times)
  \]
  induced by the commutator. Since $\bar{K}^\times$ is divisible by~$n$,
  this map is an isomorphism by a result from group cohomology. Since
  theta groups are represented on the left as the elements that map
  to the Weil paring~$e_n$ in the right hand group, there is only one
  such extension, up to isomorphism.
  
  Sketch of long, but down-to-earth proof: choosing any set-theoretic
  section $E[n] \To \Theta$ mapping $O$ to the neutral element, 
  the underlying set of~$\Theta$ can be
  identified with $\bar{K}^\times \times E[n]$. The group structure
  is then given by a map 
  \[ f : E[n] \times E[n] \To \bar{K}^\times \]
  such that 
  \[ (\lambda, T) (\lambda', T') = (\lambda \lambda' f(T, T'), T+T') \,. \]
  From the group axioms, we find that $f$ has to satisfy
  \[ f(O, T) = f(T, O) = 1\,, \qquad
     f(T, T') f(T+T', T'') = f(T, T'+T'') f(T', T'') \,.
  \]
  (I.e., $f$ is a normalized 2-cocycle.) We also have that
  \[ f(T, T') = e_n(T, T') f(T', T) \,. \]
  If we have two theta groups $\Theta_1$ and~$\Theta_2$, with multiplication 
  given by $f_1$ and $f_2$, then setting
  $f(T, T') = f_1(T, T')/f_2(T, T')$, $f$ satisfies the cocycle condition,
  and it is symmetric: $f(T, T') = f(T', T)$. We now construct a map
  $\varphi : E[n] \To \bar{K}^\times$. Set $\varphi(O) = 1$. Pick a
  basis $T, T'$ for~$E[n]$. Set
  \begin{align*}
    \varphi(T) &= \bigl(f(T,T) f(T,2T) \dots f(T,(n-1)T)\bigr)^{-1/n} \\
    \varphi(T') &= \bigl(f(T',T') f(T',2T') \dots f(T',(n-1)T')\bigr)^{-1/n}
  \end{align*} 
  with any choice of the $n$th roots (here we need that $\bar{K}^\times$
  is divisible by~$n$). Then we continue by defining
  \begin{align*}
    \varphi(mT) &= f(T,T) f(T,2T) \dots f(T,(m-1)T) \, \varphi(T)^m \\
    \varphi(mT') &= f(T',T') f(T',2T') \dots f(T',(m-1)T') \, \varphi(T')^m \\
    \varphi(mT + m'T') &= f(mT, m'T') \, \varphi(mT) \varphi(m'T')
  \end{align*}
  Now we have by an easy induction using the cocycle relation that
  \[ f(aT, bT) = \frac{\varphi(aT+bT)}{\varphi(aT)\varphi(bT)},
     f(aT', bT') = \frac{\varphi(aT'+bT')}{\varphi(aT')\varphi(bT')}\,
     f(aT, bT') = \frac{\varphi(aT+bT')}{\varphi(aT)\varphi(bT')}.
  \]
  Since we can express $f(aT+bT', cT+dT')$ in terms of values like the
  above, we get that
  \[ f(P, Q) = \frac{\varphi(P+Q)}{\varphi(P) \varphi(Q)} \]
  for all $P, Q \in E[n]$. Then 
  \[ \Theta_2 \ni (\lambda, P)
       \longmapsto (\varphi(P) \lambda, P) \in \Theta_1 
  \]
  is an isomorphism.
\end{Proof}

We deduce that $H^1(K, E[n])$ parametrizes theta groups of level~$n$ for~$E$,
up to $K$-isomorphism. These theta groups are not geometric
objects like our $n$-covering curves, but they are quite useful.

If $C \To \BP^{n-1}$ in our second interpretation represents an element
of~$\Sel^{(n)}(K, E)$, then we can easily find the corresponding theta
group~$\Theta_C$. Namely, $E[n]$ acts by automorphisms on this diagram
and thus gives us a homomorphism $\chi_C : E[n] \To \PGL_n$. As before,
we can then define $\Theta_C$ to be the pull-back of the image of~$\chi_C$
under the canonical map $\GL_n \To \PGL_n$:
\[ \SelectTips{cm}{}
   \xymatrix{ 0 \ar[r] & \Gm \ar@{=}[d] \ar[r]
                & \Theta_C \ar[d] \ar[r]
                & E[n] \ar[d]^{\chi_C} \ar[r] & 0 \\
              0 \ar[r] & \Gm \ar[r] & \GL_n \ar[r] & \PGL_n \ar[r] & 0
            }
\]
For more general diagrams $C \To S$, $\GL_n$ and~$\PGL_n$ have to be
replaced by their twists corresponding to $\Ob(C \to S) = S$; in this way,
we obtain $\Theta_C$ as a subgroup of $A_S^\times$, where $A_S$ is the
central simple algebra corresponding to~$S$. (See below.)


\section{Computation of the Selmer Group as an Abstract Group}

The interpretations given so far are not very well suited for actually
{\em computing} the $n$-Selmer group. (For $n=2$, this is not quite
true: John Cremonas {\tt mwrank} program actually enumerates 2-coverings
in order to find the 2-Selmer group.) So we will need some other
representation of the Selmer group that is more algebraic in nature
and yields itself more easily to computation.

Before we go into this, let me remark that the various Selmer groups
are related. From the diagram
\[ \SelectTips{cm}{}
   \xymatrix{ 0 \ar[r] & E[m] \ar[r] \ar[d] & E \ar[r]^m \ar@{=}[d] 
                       & E \ar[r] \ar[d]^n & 0 \\
              0 \ar[r] & E[mn] \ar[r] \ar[d]^{m} & E \ar[r]^{mn} \ar[d]^m
                       & E \ar[r] \ar@{=}[d] & 0 \\
              0 \ar[r] & E[n] \ar[r] & E \ar[r]^n & E \ar[r] & 0
            }
\]
we can deduce an exact sequence
\begin{align*}
   0 &\To \frac{E(K)[n]}{mE(K)[mn]} \To \Sel^{(m)}(K, E)
      \To \Sel^{(mn)}(K, E) \\
     &\To \Sel^{(n)}(K, E) \To \frac{\Sha(K, E)[n]}{m \Sha(K, E)[mn]} \To 0 \,.
\end{align*}
Here the map $\Sel^{(m)} \To \Sel^{(mn)}$ is induced by the inclusion
$E[m] \inj E[mn]$, and the map $\Sel^{(mn)} \To \Sel^{(n)}$ is induced
by multiplication by~$m$ $E[mn] \To E[n]$. In particular, the composition
$\Sel^{(n)} \To \Sel^{(mn)} \To \Sel^{(n)}$ (the first map coming from
the diagram where the roles of $m$ and~$n$ are exchanged) is multiplication
by~$m$ on~$\Sel^{(n)}$. Together with the exact sequence above, this
implies that
\[ \Sel^{(mn)}(K, E) \cong \Sel^{(m)}(K, E) \times \Sel^{(n)}(K, E) \]
whenever $m$ and~$n$ are coprime.

Therefore it is sufficient to compute $\Sel^{(n)}(K, E)$ when $n = p^f$
is a prime power. The first step in this is to consider the case $n = p$.
The computation of $\Sel^{(p^2)}(K, E)$ (and for higher powers of~$p$)
then is most easily done by computing the fibers of the canonical map
$\Sel^{(p^2)}(K, E) \To \Sel^{(p)}(K, E)$ one by one. (This procedure is
sometimes referred to as ``second'' or ``higher descent''.)


\subsection*{Fourth interpretation: via \'etale algebras}\strut

The goal of the interpretation I will explain now is to ease computation
of the Selmer group. So we will consider an algebraic realization 
of~$H^1(K, E[n])$ and not a geometric one.

In what follows, one of the main characters of the play will be the
{\em \'etale algebra}~$R$ of~$E[n]$. This is just the affine coordinate
algebra of the 0-dimensional scheme~$E[n]$. More concretely, we have
\[ R = \Map(E[n], \bar{K})^{G_K} \,; \]
these are $G_K$-equivariant maps on~$E[n]$ with values in~$\bar{K}$.
Note that the action of~$\sigma \in G_K$ is 
\[ f \longmapsto \bigl(f^\sigma : T \mapsto f(T^{\sigma^{-1}})^\sigma\bigr) \,. 
\]
For example, any rational function on~$E$ defined over~$K$ and not
having poles in~$E[n]$ will give an element of~$R$.

Since $E[n]$ is an \'etale $K$-scheme, $R$ is an \'etale algebra: it
is isomorphic to a product of (finite) field extensions of~$K$, one
for each $G_K$-orbit on~$E[n]$. If $T$ is a point in one such orbit,
then the corresponding field extension is $K(T)$ (i.e., $K$ with the
coordinates of the $n$-torsion point~$T$ adjoined). For example, we
always have a splitting $R = K \times R_1$, where the $K$ corresponds
to the singleton orbit~$\{O\}$, and~$R_1$ corresponds to 
$E[n] \setminus \{O\}$.

We will also use $\bar{R} = R \tensor_K \bar{K} = \Map(E[n],\bar{K})$.
As an algebra, this is just $\bar{K}^{E[n]} \cong \bar{K}^{n^2}$, but the
action of~$G_K$ is ``twisted'' by its action on~$E[n]$, permuting the
factors. 

In this context, $\bar{R}^\times$ is the multiplicative group of maps
from $E[n]$ into~$\bar{K}^\times$, and $R^\times$ is the subgroup of
$G_K$-equivariant such maps. For example, every $T \in E[n]$ gives a
map 
\[ e(T) : E[n] \ni S \mapsto e_n(T, S) \in \mu_n \subset \bar{K}^\times \,, \]
and so we obtain an injective (because $e_n$ is non-degenerate) homomorphism
\[ e : E[n] \To \bar{R}^\times \,. \]
The idea now is to extend this to (the beginning of) a resolution
of~$E[n]$ as a $K$-Galois module, in order to get some more or less
concrete realization of~$H^1(K, E[n])$. 

In general, if $R_A$ and $R_B$ are the coordinate rings of two affine
$K$-schemes $A$ and~$B$, then $R_A \tensor_K R_B$ is the coordinate
ring of $A \times B$. So $R \tensor_K R$ is the algebra of $G_K$-equivariant
maps from $E[n] \times E[n]$ into~$\bar{K}$, and 
$\bar{R} \tensor_{\bar{K}} \bar{R} = (R \tensor_K R) \tensor_K \bar{K}$
is the algebra of all such maps. 

Now observe that all $e(T)$, for $T \in E[n]$, are not only maps, but
homomorphisms $E[n] \To \bar{K}^\times$. I.e., $\alpha = e(T)$ satisfies the 
relations $\alpha(T_1 + T_2) = \alpha(T_1) \alpha(T_2)$. 
By the non-degeneracy of the Weil pairing, we also know that the $e(T)$ 
are {\em all} the homomorphisms. This implies that, if we define
\[ \partial : \bar{R}^\times \ni \alpha
        \longmapsto \Bigl((T_1, T_2)
          \mapsto \frac{\alpha(T_1) \alpha(T_2)}{\alpha(T_1+T_2)}\Bigr)
            \in (\bar{R} \tensor_{\bar{K}} \bar{R})^\times \,,
\]
then
\[ 0 \To E[n] \stackrel{e}{\To} \bar{R}^\times 
          \stackrel{\partial}{\To} (\bar{R} \tensor_{\bar{K}} \bar{R})^\times
\]
will be exact.

In order to use that to realize $H^1(K, E[n])$, we need to know what the
image of the map~$\partial$ is. One obvious property of all $\partial \alpha$
is that they are symmetric: 
$\partial\alpha(T_1, T_2) = \partial\alpha(T_2, T_1)$. But there are more
conditions they satisfy. Let us define
\[ \partial : (\bar{R} \tensor_{\bar{K}} \bar{R})^\times \ni \rho
     \longmapsto \Bigl((T_1, T_2, T_3) \mapsto
         \frac{\rho(T_1, T_2) \rho(T_1+T_2, T_3)}%
              {\rho(T_1, T_2+T_3) \rho(T_2, T_3)}\Bigr)
     \in (\bar{R} \tensor_{\bar{K}} \bar{R} \tensor_{\bar{K}} \bar{R})^\times
   \,.
\]
Then we have the following result. We define $\Sym^2_{\bar{K}}(\bar{R})$
to be the subalgebra of $\bar{R} \tensor_{\bar{K}} \bar{R}$ consisting
of symmetric maps (and similarly $\Sym^2_K(R)$ for $R$ and~$K$).

\begin{Proposition}
  The following is an exact sequence.
  \[ 0 \To E[n] \stackrel{e}{\To} \bar{R}^\times
       \stackrel{\partial}{\To} \bigl(\Sym^2_{\bar{K}}(\bar{R})\bigr)^\times
       \stackrel{\partial}{\To} 
        (\bar{R} \tensor_{\bar{K}} \bar{R} \tensor_{\bar{K}} \bar{R})^\times
  \]
\end{Proposition}

\begin{Proof}
  We only have to show exactness at 
  $\bigl(\Sym^2_{\bar{K}}(\bar{R})\bigr)^\times$. It is immediately checked
  that $\partial \partial \alpha = 1$ for $\alpha \in \bar{R}^\times$.
  On the other hand, if $\rho \in \bigl(\Sym^2_{\bar{K}}(\bar{R})\bigr)^\times$
  such that $\partial \rho = 1$, then as a map $E[n]^2 \To \bar{K}^\times$,
  $\rho$ is a symmetric 2-cocycle, and we have seen in our discussion of
  theta groups that each such 2-cocycle is a coboundary, which translates
  into $\rho = \partial\alpha$ for some~$\alpha \in \bar{R}^\times$.
\end{Proof}

\begin{Corollary} \label{RepH}
  There is an isomorphism
  \[ H = \frac{\ker\bigl(\partial \mid \Sym^2_K(R)^\times\bigr)}%
              {\partial R^\times}
         \stackrel{\cong}{\To} H^1(K, E[n]) \,.
  \]
  It is defined as follows. Take $\rho \in \Sym^2_K(R)^\times$ such
  that $\partial\rho = 1$. Then there is some $\gamma \in \bar{R}^\times$
  such that $\partial \gamma = \rho$. Now the Galois 1-cocycle
  $\sigma \mapsto \gamma^\sigma/\gamma$ takes values in the kernel 
  of~$\partial$, so we can write $\gamma^\sigma/\gamma = e(T_\sigma)$,
  where $\sigma \mapsto T_\sigma$ is a 1-cocycle with values in~$E[n]$
  representing the image of~$\rho$ in~$H^1(K, E[n])$.
\end{Corollary}

\begin{Proof}
  From the proposition, we get the short exact sequence
  \[ 0 \To E[n] \stackrel{e}{\To} \bar{R}^\times \stackrel{\partial}{\To}
       \ker\bigl(\partial \mid (\Sym^2_{\bar{K}}(\bar{R}))^\times\bigr)
       \To 0 \,.
  \]
  Apply the long exact cohomology sequence to get
  \[ R^\times \stackrel{\partial}{\To} 
       \ker\bigl(\partial \mid \Sym^2_K(R)^\times\bigr) \stackrel{\delta}{\To}
       H^1(K, E[n]) \To H^1(K, \bar{R}^\times) = 0 \,.
  \]
  (The latter equality is an easy generalization of Hilbert's Theorem~90.)
  The description of the isomorphism follows the definition of the
  connecting map~$\delta$ above.
\end{Proof}

Therefore, we can represent $H^1(K, E[n])$ as a subquotient of
$\Sym^2_K(R)^\times/(\Sym^2_K(R)^\times)^n$. Putting in the 
condition that the elements are unramified outside the set~$S$ of places
of~$K$ described in Thm.~\ref{Ramification}, we see that $\Sel^{(n)}(K, E)$
is contained in a subquotient of the ``$n$-Selmer group'' of the $S$-integers
of~$\Sym^2_K(R)$,
\[ \Sym^2_K(R)(S, n)
    = \frac{\{\rho \in \Sym^2_K(R)^\times
               : \Sym^2_K(R)(\sqrt[n]{\rho}) \text{\ unramified outside~$S$}\}}%
           {(\Sym^2_K(R)^\times)^n} \,.
\]
This is a finite group that is effectively computable. However, this
computation requires the knowledge of the class and unit groups of the
number fields occurring as factors of~$\Sym^2_K(R)$. 
If $n = p$ is an odd prime, we have a splitting
\[ \Sym^2_K(R) \cong K \times \prod_j L_j \times L \,, \]
where $K$ corresponds to $\{(O, O)\}$, $L_j$ corresponds to the set
$\{(T, jT) : O \neq T \in E[p]\}$ (where $j$ runs through a set of representatives
of~$\F_p^\times$ modulo identifying inverses), and $L$ corresponds to 
the set of unordered bases of~$E[p]$. The cardinality of this set, and
therefore the degree of~$L$, is $(p-1)^2 p (p+1)/2$. Generically, $L$
is a number field of that relative degree over~$K$, and so even for $p=3$,
this is outside the range of practical applicability of current methods
in computational algebraic number theory. So we will need a better, smaller
representation. But let us first see how one could use our current
representation at least in principle to actually compute the Selmer group.

For this, note that the result of Thm.~\ref{Ramification} can be extended
to show that the Selmer group can be obtained by a finite computation:
\[ \Sel^{(n)}(K, E)
     = \{\xi \in H^1(K, E[n]; S) :
                \forall v \in S: \res_v(\xi) \in \delta_v(E(K_v)/nE(K_v))\} 
\]
(For this, one needs to show that the image of $E(K_v)/n E(K_v)$
in~$H^1(K_v, E[n])$ is exactly the unramified part, for all $v \notin S$.)
In order to turn this into an algorithm, we first need to find the
image~$H_S$ of~$H^1(K, E[n]; S)$ in~$H$. This is obtained through
the computation of $\Sym^2_K(R)(S, n)$ and then passing to the relevant
subquotient. (For this, we need to find $R(S, n)$ and compute its
image under~$\partial$.) Then we need to have explicit representations
of the maps $\res_v$ and~$\delta_v$, for all $v \in S$. For the restriction
maps, we just observe that the result of Cor.~\ref{RepH} works over
any field~$K$ and is functorial. Applying this observation to the
field extension $K \subset K_v$, we get
\[ \res_v : H \To H_v
     = \frac{\ker\bigl(\partial \mid (\Sym^2_K(R) \tensor_K K_v)^\times\bigr)}%
            {\partial (R \tensor_K K_v)^\times} \,;
\]
this is induced by the canonical map 
$\Sym^2_K(R) \To \Sym^2_K(R) \tensor_K K_v$. Note also that $H_v$ is a
finite group. To get a nice representation of~$\delta_v$ (or~$\delta$),
we remind ourselves of the usual definition of the Weil pairing.

Let $T \in E[n]$ be an $n$-torsion point. Then there is a rational
function $G_T \in K(T)(E)^\times$ such that
\[ \diw(G_T) = n^*(T) - n^*(O)
             = \sum_{P : nP = T} P - \sum_{Q : nQ = O} Q \,.
\]
This divisor is stable under translations by elements of~$E[n]$, therefore
we have that
$G_T(P + S) = c(S) G_T(P)$ for some constant $c(S)$ independent of~$P$,
when $S \in E[n]$. This constant is $e_n(S, T)$. We can choose these
functions $G_T$ in such a way that $G : T \mapsto G_T$ is Galois-equivariant.
Then we can interpret $G$ as an element of $R(E)^\times$ (where
$R(E) = K(E) \tensor_K R$; its elements are $G_K$-equivariant maps
from $E[n]$ into $\bar{K}(E)$), and we have $G(P + T) = e(T) G(P)$
for $P \in E \setminus E[n^2]$ and $T \in E[n]$.
Now consider the following diagram.
\[ \SelectTips{cm}{}
   \xymatrix{ 0 \ar[r] & E[n] \ar[r] \ar@{=}[d]
                       & E \ar[r]^n \ar@{-->}[d]^G
                       & E \ar[r] \ar@{-->}[d]^r & 0 \\
              0 \ar[r] & E[n] \ar[r]^e 
                       & \bar{R}^\times \ar[r]^-\partial
           & \ker\bigl(\partial \mid \Sym_{\bar{K}}^2(\bar{R})^\times\bigr)
             \ar[r] & 0
            }
\]
The map $r$ is defined as shown: to find $r(P)$, take some $Q \in E$
such that $nQ = P$, then $r(P) = \partial G(Q)$. This is well defined,
since another choice $Q'$ in place of~$Q$ will differ from $Q$ by addition
of an $n$-torsion point~$T$, so 
\[ \partial G(Q') = \partial G(Q + T) = \partial (e(T) G(Q)) = \partial G(Q) \,,
\]
since $\partial e(T) = 1$. This also shows that $r$ is defined over~$K$,
so $r \in \Sym_K^2(R)(E)^\times$. To find out what function
$r_{T_1,T_2} \in \bar{K}(E)^\times$ is, let us determine its divisor.
By definition,
\[ r_{T_1,T_2}(nQ) = \frac{G_{T_1}(Q) G_{T_2}(Q)}{G_{T_1+T_2}(Q)} \,, \]
so
\[ n^*\bigl(\diw(r_{T_1,T_2})\bigr)
    = n^*(T_1) + n^*(T_2) - n^*(T_1+T_2) - n^*(O) \,.
\]
So $r_{T_1,T_2}$ has simple zeros at $T_1$ and~$T_2$ and simple poles
at~$T_1+T_2$ and~$O$: it is the function witnessing that $T_1+T_2$ is
the sum of $T_1$ and~$T_2$. So, up to normalizing constants, with respect
to a Weierstra{\ss} model of~$E$, $r_{T_1,T_2}$ is the equation of the
line joining $T_1$ and~$T_2$ divided by the equation of the (vertical)
line joining $T_1+T_2$ and~$O$. For suitable normalisation of the $G_T$,
we can take quite concretely the following. Here, the line joining
two points $T_1$ and~$T_2$, such that $T_1,T_2,T_1+T_2 \neq O$ 
(the tangent line at~$T$ 
of~$E$, when $T_1 = T_2 = T$) is supposed to have equation
\[ y = \lambda_{T_1,T_2}\,x + m_{T_1,T_2} \,. \]
Take
\[ r_{T_1,T_2} = \left\{\begin{array}{cl}
                   1 & \text{\quad if $T_1 = O$ or $T_2 = O$;} \\
                   x - x(T_1) & \text{\quad if $T_1 + T_2 = O$, $T_1 \neq O$;}\\
                   \dfrac{y - \lambda_{T_1,T_2} x - m_{T_1,T_2}}{x - x(T_1+T_2)}
                     & \text{\quad if $T_1, T_2, T_1+T_2 \neq O$.}
                        
                 \end{array}\right.
\]

Now, chasing through the definitions of the connecting homomorphisms
\[ E(K) \stackrel{\delta}{\To} H^1(K, E[n]) \text{\quad and\quad}
   H \stackrel{\cong}{\To} H^1(K, E[n]) \,,
\]
we easily find the following.

\begin{Proposition}
  The composition $E(K) \stackrel{\delta}{\To} H^1(K, E[n]) \To H$ is induced
  by $r : E(K) \setminus E[n] \To \Sym_K^2(R)^\times$.
\end{Proposition}

The ``missing'' values on $E(K)[n]$ can be obtained by a suitable
rescaling and limiting process. This result is again valid for all 
fields~$K$ and functorial in~$K$, so we can use it to compute the
local maps $\delta_v : E(K_v)/nE(K_v) \inj H_v$. Since we can easily
compute the size of the group on the left hand side --- say, $n = p$
is a prime, then
\[ \dim_{\F_p} E(K_v)/p E(K_v) = \dim_{\F_p} E(K_v)[p] + 
     \begin{cases}
       0 & \text{\quad if $v$ is finite and $v \nmid p$;} \\
       [K_v : \Q_p] & \text{\quad if $v \mid p$;} \\
       -[K_v : \R] & \text{\quad if $v$ is infinite ---}
     \end{cases}
\]
we can just pick random points in~$E(K_v)$ until their images in~$H_v$
generate a subspace of the correct size. Having found all the ``local images''
\[ J_v = \delta_v(E(K_v)/p E(K_v)) \subset H_v \,,\quad\text{for $v \in S$,} \] 
the determination
of the Selmer group as a subgroup of~$H_S$ is reduced to linear algebra
over~$\F_p$. {\em Mutatis mutandis}, this will also work for arbitrary~$n$.

Let us summarize the discussion.

\begin{Theorem}
  There is an effective procedure for computing the $n$-Selmer group
  of an elliptic curve over a number field~$K$. It requires class group
  and unit group computations in extensions of~$K$ of the form
  $K(\{T_1, T_2\})$, where $\{T_1, T_2\}$ is an unordered pair of
  $n$-torsion points of~$E$.
\end{Theorem}

As mentioned above, as it stands, this result is rather theoretical, since
the number fields that occur are too large for practical computations.
However, we can improve the situation. Restricting to $n$-torsion subgroups
in the basic exact sequence
\[ 0 \To E[n] \stackrel{e}{\To} \bar{R}^\times \stackrel{\partial}{\To}
     \ker\bigl(\partial \mid (\Sym^2_{\bar{K}}(\bar{R}))^\times\bigr)
     \To 0 \,,
\]
we obtain
\[ 0 \To E[n] \stackrel{e}{\To} \mu_n(\bar{R}) \stackrel{\partial}{\To}
     \partial \mu_n(\bar{R}) \To 0 \,.
\]
This gives us
\[ 0 \To \frac{\bigl(\partial\mu_n(\bar{R})\bigr)^{G_K}}%
              {\partial \mu_n(R)}
     \To H^1(K, E[n]) 
     \To H^1(K, \mu_n(\bar{R})) \cong R^\times/(R^\times)^n \,.
\]
The isomorphism comes in the usual way from the fact that
$H^1(K, \bar{R}^\times) = 0$.

Now we have the following nice fact (see~\cite{DSS,SchaeferStoll}).

\begin{Proposition}
  If $n = p$ is a prime, then 
  $\bigl(\partial\mu_p(\bar{R})\bigr)^{G_K} = \partial \mu_p(R)$.
\end{Proposition}

For the proof, one basically checks all possibilities for the image
of~$G_K$ in $\Aut(E[p]) \cong \GL_2(\F_p)$, the most interesting case
being when the image is a $p$-Sylow subgroup.

\begin{Remark}
  The result is not true in general for composite~$n$. For example,
  taking $n = 4$, there are 20 conjugacy classes (out of~62) of
  subgroups of~$\GL_2(\Z/4\Z)$ such that the kernel above has order
  2 or even~4. To give a very concrete example, consider the elliptic
  curve $y^2 = x^3 + x + 2/13$ over~$\Q$; then the index of
  $\partial \mu_4(R)$ in the $G_\Q$-invariants of $\partial \mu_4(\bar{R})$
  is~2. (Here the subgroup is the one problematic one of index~2.
  It occurs generically when the discriminant of the cubic is minus a
  square.)
\end{Remark}

In any case, we have a homomorphism
\[ H \stackrel{\cong}{\To} H^1(K, E[n]) \To R^\times/(R^\times)^n \,. \]
By the definition of the Kummer map 
$R^\times/(R^\times)^n \To H^1(K, \mu_n(\bar{R}))$, it is obained as follows.
Take $\rho \in \Sym_K^2(R)^\times$ representing an element of~$H$, so
$\partial \rho = 1$. Then there is $\gamma \in \bar{R}^\times$ such that
$\partial \gamma = \rho$. Now 
$\gamma^\sigma/\gamma \in \ker\partial \subset \mu_n(\bar{R})$
for all $\sigma \in G_K$, hence $\gamma^n \in R^\times$. If we change
$\rho$ into $\rho \cdot \partial \alpha$ with $\alpha \in R^\times$,
then $\gamma^n$ changes into $\gamma^n \alpha^n$, and we get a well-defined
map $H \To R^\times/(R^\times)^n$. A somewhat more explicit description
is to say that the map is induced by
\[ \Sym_K^2(R)^\times \ni \rho
   \longmapsto \bigl(T \mapsto \prod_{j=0}^{n-1} \rho(T, jT)\bigr) \in R^\times
   \,.
\]

So when $n = p$ is a prime number, then we can use the image of~$H$
in $R^\times/(R^\times)^p$ instead of $H$ itself. The advantage of this
is obvious: the field extensions of~$K$ occurring in~$R$ are usually
{\em much} smaller than the ones in~$\Sym_K^2(R)$. Generically,
$R = K \times R_1$, with $R_1$ a field extension of~$K$ of degree~$p^2-1$.
For $K = \Q$ and $p = 3$, this leads to octic number fields, where
computations are feasible.

To make this approach work, we need to know the image $\tilde{H}$
of~$H$ in~$R^\times/(R^\times)^p$, and we need to realize the map
$E[p] \stackrel{\delta}{\To} H^1(K, E[p]) \To \tilde{H}$.

The first question is answered by looking at the map $H \to \tilde{H}$:
given $\alpha \in R^\times$, if $\alpha (R^\times)^p \in \tilde{H}$, then 
$\alpha = \gamma^p$ such that $\partial\gamma \in \Sym_K^2(R)$. This
means that
\[ \tilde{H} = \frac{\{\alpha \in R^\times
                        : \partial \alpha \in (\Sym_K^2(R)^\times)^p\}}%
                    {(R^\times)^p} \,,
\]
and we can compute $\tilde{H}_S$, the image of~$H_S$, as a subgroup
of~$R(S, p)$. (There is another description of~$\tilde{H}$ that for
$p > 3$ uses smaller fields, and it only requires checking for $p$th
powers in fields of degree at most~$p^2-1$; compare~\cite{SchaeferStoll}.)

Note also that given $\alpha \in R^\times$ representing an element
of~$\tilde{H}$, we find $\rho$ representing the corresponding element
of~$H$ as $\rho = \sqrt[p]{\partial \alpha}$. By the characterization
of~$\tilde{H}$, the root exists, and since $H \cong \tilde{H}$, it does
not matter which root we take if there is a choice (as long as it is
symmetric and a cocycle); they will all represent the same element of~$H$.

For the realization of~$\delta$, consider the following diagram.
\[ \SelectTips{cm}{}
   \xymatrix{ 0 \ar[r] & E[p] \ar[r] \ar[d]^e & E \ar[r]^p \ar@{-->}[d]^G
                   & E \ar[r] \ar@{-->}[d]^F & 0 \\
              0 \ar[r] & \mu_p(\bar{R}) \ar[r] & \bar{R}^\times \ar[r]^p
                   & \bar{R}^\times \ar[r] & 0
            }
\]
Here, $F \in R(E)^\times$ is the function such that $F_T(pQ) = G_T(Q)^p$ 
for $T \in E[p]$.
We find that the divisor of~$F_T$ is $p \cdot T - p \cdot O$, and if
$F(pQ) = G(Q)^p$ with $G \in R(E)^\times$, then $F$ induces a well-defined
map $F : E(K) \setminus E[p] \To R^\times/(R^\times)^p$, independent of
the particular choice of~$F$. Tracing through the definitions shows:

\begin{Proposition}
  The composition 
  $E(K) \stackrel{\delta}{\To} H^1(K, E[p]) \To R^\times/(R^\times)^p$
  is given by~$F$ on $E(K) \setminus E[p]$.
\end{Proposition}

Again, this is functorial in~$K$, and so we can use it for the local
maps~$\delta_v$. The algorithm for computing $\Sel^{(p)}(K, E)$ then
works as before, but now working within~$R$ instead of~$\Sym_K^2(R)$.

\begin{Theorem}
  There is an algorithm that computes $\Sel^{(p)}(K, E)$, which is
  efficient modulo computation of class and unit groups in the number
  fields $K(T)$, where $T$ runs through points of order~$p$ on~$E$.
\end{Theorem}


\section{Constructing geometric representations \\ of Selmer group elements}

Our goal in the following will be to find explicitly the $n$-coverings 
corresponding to given elements of the $n$-Selmer group. We assume that
we have realized the Selmer group as a subgroup of
\[ H = \frac{\ker\bigl(\partial \mid \Sym_K^2(R)^\times\bigr)}%
            {\partial R^\times} \,.
\]
In practice, $n = p$, and we will have computed the $p$-Selmer group
as a subgroup of~$\tilde{H} \subset R^\times/(R^\times)^p$, but we
can easily transfer this to~$H$, by the map 
$\alpha \mapsto \sqrt[p]{\partial \alpha}$ on representatives.

We first need to explain the connection between our third and fourth
interpretations of elements of~$H^1(K, E[n])$ in some detail.
By the general theory of central group extensions, theta groups
are classified by the $G_K$-equivariant symmetric 2-cocycles 
in~$Z^2(E[n], \bar{K}^\times)$, modulo the coboundaries of 
$G_K$-equivariant 1-cochains. The correspondence
is as follows. First note that for every theta group
\[ 0 \To \Gm \stackrel{\alpha}{\To} \Theta \To E[n] \To 0 \,, \]
there is a $K$-defined set-theoretic section $E[n] \To \Theta$.
To see this, pick any section $s : E[n] \To \Theta$; then for any
$\sigma \in G_K$, $s^\sigma/s$ gives a map $E[n] \To \Gm$, which (as
a function of~$\sigma$) is
a cocycle, hence can be interpreted as an element of~$Z^1(K, \bar{R}^\times)$.
Now $H^1(K, \bar{R}^\times) = 0$, so there is some map $t : E[n] \To \Gm$
such that $s^\sigma(T) t^\sigma(T) = s(T) t(T)$; replacing $s$ by~$st$
therefore yields a $K$-defined section.

Given such a section~$s$, we obtain a $K$-defined 2-cocycle 
$\phi : E[n]^2 \To \Gm$ by setting 
$\phi(T_1, T_2) = \alpha^{-1}(s(T_1) s(T_2) s(T_1+T_2)^{-1})$. 
Changing the section~$s$ amounts to a change of~$\phi$ by the coboundary
of a $K$-defined 1-cochain. The commutator condition translates into
$\phi(T_1, T_2) = e_n(T_1, T_2) \phi(T_2, T_1)$.

Now if we fix a $K$-defined section 
$\tilde{\chi}_E : E[n] \To \Theta_E \subset \GL_n$, then we obtain
a 2-cocycle $\varepsilon$ in the way described for a general theta group.
Then the ``difference'' $\phi/\varepsilon$ will be a {\em symmetric}
2-cocycle, since the commutator condition cancels.

If $n$ is odd, then there is a specific way of choosing a lift
$\tilde{\chi}_E$ such that $\varepsilon$ becomes a power of the
Weil pairing~$e_n$; in fact, $\varepsilon = e_n^k$ such that 
$2k \equiv 1 \bmod n$.

Given a symmetric 2-cocycle $\rho$, we get from $\Theta_E$ to~$\Theta_\rho$
by ``twisting'' the multiplication in~$\Theta_E$ by~$\rho$.
Writing $(T, \lambda)$ for the element $\lambda \tilde{\chi}_E(T)$,
the multiplication in~$\Theta_E$ is
\[ (T_1, \lambda_1) (T_2, \lambda_2)
     = (T_1 + T_2, \lambda_1 \lambda_2\,\varepsilon(T_1, T_2)) \,,
\]
whereas the multiplication in~$\Theta_\rho$ will be
\[ (T_1, \lambda_1) (T_2, \lambda_2)
     = (T_1 + T_2, \lambda_1 \lambda_2\,\rho(T_1, T_2) \varepsilon(T_1, T_2)) 
    \,.
\]

Now, our definition of the maps $\partial$ was exactly such that they
correspond to the coboundary maps in the standard cochain complex
\[ C^1(E[n], \bar{K}^\times) \stackrel{\partial}{\To} C^2(E[n], \bar{K}^\times)
    \stackrel{\partial}{\To} C^3(E[n], \bar{K}^\times) \,.
\]
Therefore, the elements of~$H$ represent exactly the $K$-defined symmetric
2-cocycles modulo the coboundaries of $K$-defined 1-cochains. Tracing
through the definitions, we see that the theta group $\Theta_\rho$
corresponds to the same element of~$H^1(K, E[n])$ as the image of~$\rho$
in~$H$.

Recall the obstruction map
\[ \Ob : H^1(K, E[n]) \To H^1(K, \PGL_n) \cong H^2(K, \mu_n) = \Br(K)[n] \,. \]
In our second interpretation, it was given by mapping $C \To S$ to the
element of~$H^1(K, \PGL_n)$ corresponding to the twist~$S$ of~$\BP^{n-1}$.
From this it is obvious that $\Ob = \chi_{E,*}$ is the map induced
by $\chi_E : E[n] \To \PGL_n$ on cohomology.

Now $H^1(K, \PGL_n)$ also classifies $K$-isomorphism classes of
{\em central simple algebras} of dimension~$n^2$ over~$K$; these are
twists of the matrix algebra~$\Mat_n(K)$. Therefore, one possible way
of representing the obstruction explicitly is through the construction
of the corresponding central simple algebra. Given a theta group~$\Theta$,
we obtain this in a very simple way: observe that the set $\bar{A}_\Theta$ 
of all linear
combinations of elements of~$\Theta$ (where we use the scalar multiplication
coming from the theta group structure) is in a natural way a $\bar{K}$-algebra
of dimension~$n^2$ carrying an action of~$G_K$ (we simply extend the
multiplication we have on~$\Theta$ linearly). We let $A_\Theta$ denote
the $K$-algebra of $G_K$-invariant elements. For $\Theta_E$, we obtain
in this way the matrix algebra $\Mat_n(\bar{K})$ with its usual $G_K$-action;
this is because $\Theta_E$ naturally sits inside $\GL_n$ and spans the
matrix algebra. The $\bar{K}$-isomorphism between $\Theta_E$ and~$\Theta$
extends to a $\bar{K}$-isomorphism between $\Mat_n(\bar{K})$ 
and~$\bar{A}_\Theta$, showing that $A_\Theta$ is indeed a central simple
$K$-algebra. It is then obvious that $A_\Theta$ corresponds to
$\Ob(\Theta)$.

As an aside, note that there is a completely natural and coordinate-free
$K$-defined isomorphism
\[ \langle \Theta_E \rangle \stackrel{\cong}{\To} \End(L(n \cdot O)) \]
given by identifying $\Theta_E$ with the set of pairs $(T, f_T)$ as before and
using their action on~$L(n \cdot O)$.

How do we realize the central simple algebra $\Ob(\xi)$ in terms of
$\rho \in \Sym_K^2(R)^\times$ representing $\xi \in H^1(K, E[n])$?
Note that once we fix a section $s : E[n] \To \Theta$, the elements
of $A_\Theta = \langle \Theta \rangle^{G_K}$ are identified with
$K$-equivariant maps $E[n] \To \bar{K}$ and therefore can be viewed
as elements of~$R$. The map is
\[ A_\Theta \ni \sum_T z(T) s(T) \longmapsto (z : T  \mapsto z(T)) \in R \,. \]
Therefore, we can use $R$ as the underlying $K$-vector space, and
we only have to define a new multiplication. Let us do it first with 
$\Theta_E$. From the definitions, we get that the multiplication 
on $A = A_{\Theta_E}$ must be defined by
\[ z_1 \ast_\varepsilon z_2
    = \bigl(T \mapsto
            \sum_{T_1+T_2 = T} \varepsilon(T_1,T_2) z_1(T_1) z_2(T_2)\bigr) \,.
\]
In general, for $A_\rho = A_{\Theta_\rho}$, we define the multiplication as
\[ z_1 \ast_{\varepsilon\rho} z_2
    = \bigl(T \mapsto
            \sum_{T_1+T_2 = T} \varepsilon(T_1,T_2) \rho(T_1,T_2)
                                z_1(T_1) z_2(T_2)\bigr) \,.
\]

\begin{Proposition}
  Let $\rho = \partial \gamma$ with $\gamma \in \bar{R}^\times$. Then
  in the realizations given above, a $\bar{K}$-isomorphism between
  $A_\rho$ and~$A$ is given by 
  \[ \phi_\gamma : \bar{A}_\rho \ni z \longmapsto \gamma\,z \in \bar{A} \,, \]
  where the multiplication is that of~$\bar{R}$(!).
\end{Proposition}

\begin{Proof}
  We compute
  \begin{align*}
    \phi_\gamma(z_1 *_{\varepsilon\rho} z_2)
      &= \bigl(T \mapsto \gamma(T)
               \sum_{T_1+T_2 = T} \varepsilon(T_1,T_2) \rho(T_1,T_2)
                                    z_1(T_1) z_2(T_2)\bigr) \\
      &= \bigl(T \mapsto
               \sum_{T_1+T_2 = T} \varepsilon(T_1,T_2) \gamma(T_1) z_1(T_1)
                                   \gamma(T_2) z_2(T_2)\bigr) \\
      &= \phi_\gamma(z_1) *_\varepsilon \phi_\gamma(z_2) \,.
  \end{align*}
\end{Proof}

\begin{Remark}
  One can view $A$ and $A_\rho$ as twisted versions of the group algebra
  of~$E[n]$. The group algebra (over~$\bar{K}$) is $\bar{R}$, with 
  convolution as multiplication:
  \[ z_1 * z_2 = \bigl(T \mapsto \sum_{T_1+T_2 = T} z_1(T_1) z_2(T_2)\bigr) \,.
  \]
  We can also consider the $G_K$-invariant subalgebra $(R, *)$. Now it turns
  out that $(R, *)$ is actually isomorphic to $(R, \cdot)$ (i.e., with
  point-wise multiplication). This isomorphism is given by ``Fourier transform''
  in the following way. Define
  \[ \bar{R} \ni \alpha \longmapsto 
   \Bigl(\hat{\alpha} : T \mapsto \frac{1}{n^2} \sum_S e_n(T,S) \alpha(S)\Bigr)
      \in \bar{R} \,.
  \]
  Then one can easily check that 
  $\widehat{\alpha \beta} = \hat{\alpha} * \hat{\beta}$ and that
  $\hat{\hat{\alpha}} = \alpha/n^2$. Furthermore, $\hat{\cdot\mathstrut}$
  is defined over~$K$, so it gives an isomorphism 
  $(R, \cdot) \stackrel{\cong}{\To} (R, *)$.
  
  We get $A = (R, *_\varepsilon)$ by twisting convolution in such a way 
  that commutators on the image of $E[n]$ evaluate to the Weil pairing.
\end{Remark}

Now suppose that $\rho$ represents an element~$\xi$ in the $n$-Selmer group.
Then we know that the obstruction vanishes, hence there is an
isomorphism $\iota : A_\rho \stackrel{\cong}{\To} \Mat_n(K)$. 
Given such an isomorphism, we can write the cocycle representing
$\Ob(\xi) \in H^1(K, \PGL_n)$ explicitly as a coboundary. In the
diagram below, we obtain an automorphism (over~$\bar{K}$) of~$\Mat_n$,
which must be conjugation by some matrix $M \in \GL_n(\bar{K})$.
\[ \SelectTips{cm}{}
   \xymatrix{ \bar{A}_\rho \ar[r]^-\iota \ar[d]^{\phi_\gamma}
                & \Mat_n \ar[d]^M \\
              \bar{A} \ar[r]^-\cong & \Mat_n 
            }
\]
We can find~$M$ (which is well-defined up to a multiplicative constant)
from the automorphism by linear algebra. 

\begin{Proposition}
  We have that for all $\sigma \in G_K$,
  \[ \BP(M^\sigma M^{-1}) = \chi_E(\xi_\sigma) \,. \]
  Here $\BP(X)$ denotes the image of~$X \in \GL_n$ in~$\PGL_n$.
\end{Proposition}

\begin{Proof}
  The map labeled~$M$ in the diagram above is $X \mapsto M X M^{-1}$.
  Applying $\sigma$ to the diagram, we see that on the one hand,
  by the bottom isomorphism,
  \[ z \frac{\gamma^\sigma}{\gamma}
      \longmapsto \sum_T z(T) \frac{\gamma^\sigma}{\gamma}(T) \tilde{\chi}_E(T)
      = \sum_T z(T) e_n(\xi_\sigma, T) \tilde{\chi}_E(T) \,.
  \]
  (Recall that in our fourth interpretation, 
  $e(\xi_\sigma) = \gamma^\sigma/\gamma$, where $\partial \gamma = \rho$.)
  On the other hand, we need to have that
  \[ z \frac{\gamma^\sigma}{\gamma}
      \longmapsto \sum_T z(T) M^\sigma M^{-1} \tilde{\chi}_E(T) M M^{-\sigma}
      \,.
  \]
  Comparing these expressions shows that
  \[ M^\sigma M^{-1} \tilde{\chi}_E(T)
      = e_n(\xi_\sigma, T) \tilde{\chi}_E(T) M^\sigma M^{-1} 
  \]
  for all $T \in E[n]$. If we write 
  $M^\sigma M^{-1} = X \tilde{\chi}_E(\xi_\sigma)$, then we find that
  $X$ commutes with all $\tilde{\chi}_E(T)$ and therefore with all
  of~$\Mat_n$. Therefore $X$ must be a scalar matrix, and
  \[ \BP(M^\sigma M^{-1}) = \BP(\tilde{\chi}_E(\xi_\sigma)) 
                          = \chi_E(\xi_\sigma) \,.
  \]
\end{Proof}

From this, we can draw the following useful conclusion.

\begin{Corollary}
  With the notations above, the second interpretation 
  of an element~$\xi \in \Sel^{(n)}(K, E)$ can be realized in the form
  \[ \SelectTips{cm}{}
     \xymatrix{ C \ar[r] \ar@{-->}[d] & \BP^{n-1} \ar@{-->}[d]^M \\
                E \ar[r]^-{|n \cdot O|} & \BP^{n-1}
              }
  \]
  Note that for $n \ge 3$, the horizontal arrows are embeddings, and
  so the map on the left is given by restriction of the map on the right.
\end{Corollary}

\begin{Proof}
  We note that by the previous proposition, the cocycle associated to
  the diagram is exactly~$\xi$.
\end{Proof}

In practical terms, this means that we make the linear change of variables
corresponding to~$M^\top$ (acting on the right) in the equations of 
$E \subset \BP^{n-1}$ to obtain equations for~$C$. 

However, recall that in order to find~$M$, we need to actually {\em have}
an isomorphism between $A_\rho$ and~$\Mat_n(K)$. To find such an
isomorphism explicitly is a nontrivial problem, even though we already
know that such an isomorphism exists. When $n = 2$, this turns out to
be equivalent to finding a $K$-rational point on a conic, knowing that
it has points everywhere locally. The problem of ``trivializing the algebra''
that comes up here can be viewed as a generalization
of this very classical problem.

At least in theory, the problem can be reduced to solving a norm equation
(over~$K$ or some extension of~$K$). In practice, however, the data
defining the algebra structure on~$A_\rho$ can be very large (they come
in the end from elements of $R(S, n)$ and often will involve units
of the number fields occurring in~$R$), making this approach impractical.
On the other hand, at least when working over~$\Q$, we have a method
that seems to work very well in practice (at least when $n = 3$, which
is the only case where the first step, the computation of the Selmer group
as a group, can be carried out successfully), although we have so far
not proved that it really always works. The idea is to first find a
maximal order in~$A_\rho$, which we know is isomorphic to~$\Mat_n(\Z)$.
Then we apply a certain reduction procedure with the goal of reducing
the structure constants in size until they are small enough to read off
an isomorphism.

\medskip

Even though it is not immediately relevant to what we are doing in these
lectures, I would like to mention the following.

\begin{Proposition}\strut
  \begin{enumerate}
    \item
      The obstruction map is even: $\Ob(-\xi) = \Ob(\xi)$.
    \item
      The Weil pairing cup-product pairing
      \[ \cup_e : H^1(K, E[n]) \times H^1(K, E[n]) \To H^2(K, \mu_n) \]
      can be expressed in terms of the obstruction map:
      \[ \xi \cup_e \eta = \Ob(\xi+\eta) - \Ob(\xi) - \Ob(\eta) \,. \]
  \end{enumerate}
  In particular, $\Ob$ is a quadratic map:
  \[ \Ob(m \xi) = m^2 \Ob(\xi) \,, \qquad
     \Ob(\xi+\eta) + \Ob(\xi-\eta) = 2 \Ob(\xi) + 2 \Ob(\eta) \,.
  \]
\end{Proposition}

\begin{Proof}
  (1) We get the diagram corresponding to~$-\xi$ from the diagram
      corresponding to~$\xi$ by composing it with inversion on~$E$:
      \[ \SelectTips{cm}{}
         \xymatrix{ C \ar@{-->}[d] \ar[r] & S \ar@{-->}[d] \\
                    E \ar[d]^{-1} \ar[r] & \BP^{n-1} \ar[d] \\
                    E \ar[r] & \BP^{n-1}
                  }
      \]
      In particular, the Brauer-Severi variety for~$-\xi$ is the same
      (namely~$S$) as the one for~$\xi$. \\
  (2) This is a straight-forward, if somewhat tedious, verification
      using cocycles.
\end{Proof}

\begin{Remark}
  There is another way of obtaining the model of~$C$ in~$\BP^{n-1}$
  (which is the one actually currently used in my 3-descent prorgam).
  It roughly works as follows.
  
  Instead of mapping $E \stackrel{\phi}{\To} \BP^{n-1}$, we can also map 
  to the dual curve, i.e., we send $P \in E$ to the point 
  in~$(\BP^{n-1})^\vee$ corresponding to the osculating hyperplane
  at~$\phi(P)$ (for $n=3$, this is just the tangent line, for $n=2$,
  it is $\phi(P) \in (\BP^1)^\vee = \BP^1$ itself). We can combine $\phi$
  and this morphism $\phi^\vee : E \To (\BP^{n-1})^\vee$ into a single
  morphism and then follow it by the Segre embedding, where we can 
  identify $\BP^{n^2-1}$ with $\BP(\Mat_n)$ and the embedding with 
  multiplication of column vectors by row vectors. The image of Segre
  is therefore the set of rank-1 matrices.
  \[ E \stackrel{(\phi,\phi^\vee)}{\To} \BP^{n-1} \times (\BP^{n-1})^\vee
       \stackrel{\text{Segre}}{\To} \BP^{n^2-1} = \BP(\Mat_n) 
  \]
  The image of~$E$ in~$\BP(\Mat_n)$ will also be contained in the hyperplane
  of trace-zero matrices; this corresponds to the fact that $\phi(P)$
  is on the hyperplane~$\phi^\vee(P)$.
  
  Now there is a commutative diagram
  \[ \SelectTips{cm}{}
     \xymatrix{ E \ar[r]^-{(\phi,\phi^\vee)} \ar[d]^{G} &
                   \BP^{n-1} \times (\BP^{n-1})^\vee \ar[r]^-{\text{Segre}} &
                   \BP(\Mat_n) \\
                \BP(\bar{R}) \ar@{-->}[rr]^{\cdot t} & &
                  \BP(\bar{R}) \ar[u]^{\tilde{\chi}_E}
              }
  \]
  where $t \in R$ is a certain element satisfying $t(O) = 0$,
  $t(T) \neq 0$ for all $T \neq O$. For example, when $n = 3$ and
  $\tilde{\chi}_E$ is chosen so as to make $\varepsilon = e_3^2$,
  then $t(T) = 1/y(T)$ (w.r.t.\ a Weierstra{\ss} model of~$E$).
  
  Twisting by the cocycle represented by $\rho = \partial\gamma$,
  we obtain a model $C_1$ of $C$ in $\BP(\bar{R})$ as $\gamma^{-1} \cdot G(E)$.
  This can be computed explicitly in terms of~$\rho$ only: applying~$\partial$,
  we have 
  \[ z \in C_1 \iff \gamma z \in G(E) \iff \rho \partial z \in r(E) \,, \]
  and the latter leads to quadratic equations in~$z$, from which $r(E)$
  can be eliminated.
  
  Then $t \cdot C_1$ will be contained in the rank-1, trace-0 locus
  of $\BP(\bar{A}_\rho)$ (identifying underlying vector spaces). If we
  have an isomorphism $A_\rho \stackrel{\iota}{\To} \Mat_n(K)$, then
  we obtain $C$ by projecting $\iota(t \cdot C_1) \subset \BP(\Mat_n)$
  to any nonzero column.
\end{Remark}


\section{Minimization and Reduction}

I would like to come back to the diagram 
\[ \SelectTips{cm}{}
   \xymatrix{ \BP^1 & E \ar@{=}[d] \ar[l]_x
                           & C \ar[l]_\pi \ar@{-->}[d]^\cong \ar[r]^{|D|}
                           & S \ar[r]^-\iota_-{\cong} \ar@{-->}[d]^\cong
                           & \BP^{n-1} \\
               & E & E \ar[l]_{n} \ar[r]^-{|n \cdot O|}
                & \BP^{n-1} &
            }
\]
From this diagram, one can read off that $\pi^*(O) \sim n D$ as divisors
on~$C$. Now there is a theory of heights on varieties. For a very ample
divisor $D$, one defines 
\[ h_D : C \stackrel{|D|}{\To} \BP^N \stackrel{h}{\To} \R_{\ge 0} \]
via the logarithmic height on~$\BP^N$. This is well-defined up to bounded
functions. The main facts are that this induces a homomorphism
\[ \{\text{divisors on~$C$}\}
        \To \frac{\{\text{functions $C \to \R_{\ge 0}$}\}}%
                                       {\{\text{bounded functions}\}}
\]
and that it is compatible with dominant morphisms $C' \stackrel{\phi}{\To} C$ 
in the sense that
\[ h_{\phi^* D} = h_D \circ \phi + O(1) \]
for divisors $D$ on~$C$.

So if $Q \in C(K)$ and $P = \pi(Q) \in E(K)$, we get that
\begin{align*}
  h_D(Q) &= \frac{1}{n} h_{\pi^* O}(Q) + O(1)
          = \frac{1}{n} h_O(P) + O(1) \\
         &= \frac{1}{2n} h_{2O}(P) + O(1)
          = \frac{1}{2n} h(x(P)) + O(1) \,.
\end{align*}
So up to something bounded, going from $E$ to~$C$ divides logarithmic
heights by~$2n$.

Now, to really make use of this, one needs the ``$O(1)$'' to be fairly
small. How large the error is depends mainly on the size of the coefficients
in the equations describing $C$ (and~$E$) --- the smaller they are,
the better. So we would like to choose coordinates on~$\BP^{n-1}$ in
such a way that $C$ is described by small equations.

There are two ways in which the equations can be large. The first is that
the model is not minimal, i.e., it has unnecessary prime powers in its
discriminant. So in a first step, one will try to remove these and obtain
a minimal model. This step, called ``minimization'', has been worked out in 
theory and practice for $n = 2,3,4$ and is being worked on for~$n=5$.

Then, assuming now that the model is minimal, we still can make coordinate
changes by $\SL_n(\Z)$. So we would like to find such a coordinate change
that makes the coefficients of our equations small. This step, called
``reducion'', has been worked
out (for this special case at least) in theory for all~$n$, and is 
implemented for $n = 2,3,4$.

\medskip

There is an implementation in {\sf MAGMA} that computes the 3-Selmer group
$\Sel^{(3)}(\Q, E)$ (inside $R^\times/(R^\times)^3$) of an elliptic curve
$E$ over~$\Q$. It then transfers the elements into~$H$, finds the structure
constants for the central simple algebras associated to them, computes
an isomorphism with $\Mat_3(\Q)$, finds the equations for the model
in~$\BP^8$ and then projects the model into~$\BP^2$; finally this is
minimized and reduced. This program works quite well in practice for
curves if moderate size and produces a list of curves corresponding to
$(\Sel^{(3)}(\Q, E) \setminus \{0\})/\{\pm 1\}$. (Note that elements that
are negatives of each other give rise to the same curve, which has two
different structures as a principal homogeneous space for~$E$.)



\begin{thebibliography}{MMMMM9}
\frenchspacing

\bibitem[CF]{CasselsFroehlich}
  {\sc J.W.S.~Cassels, A.~Fr\"ohlich:} {\it Algebraic number theory,}
  Academic Press (1993).
\bibitem[CFOSS1]{Paper1}
  {\sc J.E.~Cremona, T.A.~Fisher, C.~O'Neil, D.~Simon, M.~Stoll:}
  {\it Explicit $n$-descent on elliptic curves. I. Algebra,}
  Preprint (2006).
\bibitem[CFOSS2]{Paper2}
  {\sc J.E.~Cremona, T.A.~Fisher, C.~O'Neil, D.~Simon, M.~Stoll:}
  {\it Explicit $n$-descent on elliptic curves. II. Geometry,}
  Preprint (2006).
\bibitem[CFOSS3]{Paper3}
  {\sc J.E.~Cremona, T.A.~Fisher, C.~O'Neil, D.~Simon, M.~Stoll:}
  {\it Explicit $n$-descent on elliptic curves. III. Algorithms,}
  in preparation.
\bibitem[DSS]{DSS}
  {\sc Z.~Djabri, E.F.~Schaefer, N.~Smart:} {\it Computing the
  $p$-Selmer group of an elliptic curve}, 
  Trans.\ Amer.\ Math.\ Soc.\ {\bf 352} (2000), 5583--5597.
\bibitem[ScSt]{SchaeferStoll}
  {\sc E.F.~Schaefer, M.~Stoll:} {\it How to do a $p$-descent on an elliptic 
   curve}, Trans. Amer. Math. Soc. {\bf 356}, 1209--1231 (2004).
\bibitem[Sil]{Silverman}
  {\sc J.H. Silverman:} {\it The arithmetic of elliptic curves},
  Springer GTM 106 (1986).

\end{thebibliography}
\end{document}